\documentclass[a4paper,10pt]{amsart}
\usepackage{amsmath,amssymb,amsthm}

\newcommand{\CM}{\overline{\operatorname{\mathcal{M}}}}
\newcommand{\CB}{\overline{\operatorname{\mathcal{B}}}}
\newcommand{\CE}{\overline{\operatorname{\mathcal{E}}}}

\newcommand{\IM}{\operatorname{\mathcal{M}}}
\newcommand{\BB}{\operatorname{\mathcal{B}}}
\newcommand{\EE}{\operatorname{\mathcal{E}}}

\newcommand{\Si}{\dot{S}}

\newcommand{\CR}{\bar{\partial}}

\newcommand{\Aut}{\operatorname{Aut}}

\newcommand{\Crit}{\operatorname{Crit}}

\newcommand{\loc}{\operatorname{loc}}

\newcommand{\ind}{\operatorname{ind}}
\newcommand{\coker}{\operatorname{coker}}
\newcommand{\del}{\partial}

\newcommand{\Ju}{\underline{J}}

\newcommand{\IC}{\operatorname{\mathbb{C}}}

\newcommand{\IR}{\operatorname{\mathbb{R}}}

\title{Transversality problems in \\symplectic field theory \\and a new Fredholm theory}
\author{Oliver Fabert}

\begin{document}

\maketitle

\begin{abstract}
This survey wants to give a short introduction to the transversality problem in symplectic field theory and motivate to approach 
it using the new Fredholm theory by Hofer, Wysocki and Zehnder. With this it should serve as a lead-in to the 
user's guide to polyfolds, which is appearing soon and the result of a working group organized by J. Fish, 
R. Golovko and the author at MSRI Berkeley in fall 2009. \\

\end{abstract}

\tableofcontents

\markboth{O. Fabert}{Transversality problems in SFT} 

\section*{Summary}

Many problems in symplectic topology, like Gromov's pseudoholomorphic curves and Floer homology, are based on the study 
of solution sets of nonlinear elliptic PDE. These moduli spaces exhibit lack of compactness, but on the other hand lead 
to very elaborate compactifications, which are the source for interesting algebraic invariants. In many cases, including 
symplectic field theory (SFT), the algebraic structures of interest are precisely those created by the violent analytic 
behavior. \\

In order to provide the right framework in order to deal with all these problems, H. Hofer, K. Wysocki and E. Zehnder 
have introduced a new class of spaces, called polyfolds, which serve as new ambient space containing the compactified 
moduli space as the zero set of the nonlinear Cauchy-Riemann operator, for which they provide the usual Fredholm package 
including an implicit function theorem and abstract perturbation scheme, see [HWZ] and [H] and the references therein. \\

As a result of the working group on polyfolds which the author organized together with Joel Fish and Roman Golovko 
at the MSRI in Berkeley in fall 2009, the three organizers use their own lecture notes in order to produce a 
user's guide to polyfolds, whose main goal is to equip symplectic geometers (and others) with a relatively short 
guide to show them how to apply the polyfold theory to their own Fredholm problem. \\

This survey paper wants to give a short introduction to the transversality problem in SFT and hence serve as a lead-in 
for the user's guide which will appear soon. \\

After a short introduction 
to symplectic field theory in section 1, which will serve as our main example in this introduction, we will first recall 
in section 2 the classical 
approach to equip moduli spaces and their compactifications with nice manifold structures (with boundaries and corners) 
using infinite-dimensional Banach space bundles over Banach manifolds. Since everything relies on an infinite-dimensional 
version of the classical implicit function theorem, it turns out that the crucial step is to prove a transversality result 
for the Cauchy-Riemann operator. While it is well-known that transversality holds for a generic choice of almost complex 
structure as long as all holomorphic curves are simple, that is, not multiply-covered, we will discuss in section 3 how 
severe the transversality problem with multiply-covered curves actually is, where it will turn out that the biggest 
problems arise 
from branched covers of orbit cylinders, also called orbit curves. Instead of discussing other approaches like Kuranishi 
structures and virtual fundamental cycle techniqes which claim to prove the transversality problem for holomorphic curves 
in full generality, we will shortly discuss in section 4 special approaches 
to the transversality problem like automatic transversality, obstruction bundles and domain-dependent almost complex 
structures, which are applicable in special cases and turn out to be useful when one is interested to perform actual 
computations. After seeing the limitations in all existing approaches (except the ones not discussed), we hope that 
the reader is finally motivated enough to approach the transversality problem using a completely new Fredholm theory 
based on new geometric objects like sc-manifolds and polyfolds, where will outline the main ideas and new concepts in the 
last section 5. \\  

Since the goal of this survey is just to show how severe the transversality problem in SFT actually is and thereby motivate symplectic 
geometers to get in touch with the results of the great polyfold project, the author wants to apologize right away for all the appearing 
simplifications, inaccuracies and missing literature. On the other hand, comments, remarks and suggestions are very 
welcome. \\
    
The author would like to thank H. Hofer, K. Wehrheim, his co-organizers J. Fish and R. Golovko and the audience of 
their working group at the MSRI for their great help, stimulating discussions and useful comments. Since this survey 
does not claim to provide any original new results, he further wants to thank K. Wysocki and E. Zehnder but also 
C. Wendl, K. Cieliebak, K. Mohnke, E. Ionel, C. Taubes and others for their great work towards transversality for 
holomorphic curves. Since this survey was written up when the author was a postdoc at the Max Planck Institute (MPI)
for Mathematics in the Sciences, he wants to thank both the MSRI and the MPI for the hospitality and their great working environment.      

\section{A short introduction to symplectic field theory (SFT)}
Symplectic field theory (SFT), introduced by H. Hofer, A. Givental and Y. Eliashberg in 2000 ([EGH]), is a very large 
project and can be viewed as a topological quantum field theory approach to Gromov-Witten theory. Besides providing a 
unified view on established pseudoholomorphic curve theories like symplectic Floer homology, contact homology and 
Gromov-Witten theory, it leads to numerous new applications and opens new routes yet to be explored. \\

In this survey we will restrict our attention to the transversality problem in SFT. Apart from the fact that 
SFT aims at being a grand unified theory of $J$-holomorphic curves in symplectic and contact topology, the main 
motivation comes from the observation that, in contrast to Gromov-Witten theory and symplectic Floer homology, there 
do not exist reasonable assumptions on the target manifolds like monotonicity or, more generally, semipositivity, for 
which transversality (for all relevant moduli spaces) can be proven for generic choices of almost complex structures. 
Furthermore, while in Gromov-Witten 
theory K. Cieliebak and K. Mohnke were able to prove in [CM] a transversality result for general symplectic manifolds 
using domain-dependent almost complex structures and employing Donaldson's construction of symplectic hypersurfaces, a 
corresponding result in SFT was so far only established in the so-called Floer case by the author in [F1] and 
currently seems out of reach for the general case, as speculated by K. Cieliebak in his talk in the "broken dreams" seminar 
at MSRI in fall 2009.    \\

We start with briefly recalling the geometric setup of SFT for closed contact manifolds $(V,\xi=\{\lambda=0\})$. 
For the moduli spaces and the functional analytic setup we will for simplicity only consider the case of genus zero. \\
   
Recall that a contact one-form $\lambda$ defines a vector field $R$ on $V$ by 
$R\in\ker d\lambda$ and $\lambda(R)=1$, which 
is called the Reeb vector field. We assume that 
the contact form is Morse in the sense that all closed orbits of the 
Reeb vector field are nondegenerate in the sense of [BEHWZ]; in particular, the set 
of closed Reeb orbits is discrete. \\

The SFT invariants are defined by counting 
$\Ju$-holomorphic curves in $\IR\times V$ which are asymptotically cylindrical over 
chosen collections of Reeb orbits $\Gamma^{\pm}=\{\gamma^{\pm}_1,...,
\gamma^{\pm}_{n^{\pm}}\}$ as the $\IR$-factor tends to $\pm\infty$, see [BEHWZ]. 
The almost complex structure $\Ju$ on the cylindrical 
manifold $\IR\times V$ is required to be cylindrical in the sense that it is  
$\IR$-independent, links the two natural vector fields on $\IR\times V$, namely the 
Reeb vector field $R$ and the $\IR$-direction $\del_s$, by $\Ju\del_s=R$, and turns 
the distribution $\xi$ on $V$ into a complex subbundle of $TV$, 
$\xi=TV\cap \Ju TV$. \\

Then the moduli space $\IM^0=\IM^0(\Gamma^+,\Gamma^-)$ of parametrized curves consists of tuples 
$(u,j)$ where $j$ is a complex structure on the punctured Riemann sphere $\Si$ with $s^+=\#\Gamma^+$ positive 
and $s^-=\#\Gamma^-$ negative punctures and $u:(\Si,j)\to(\IR\times V,\Ju)$ is a $\Ju$-holomorphic map in the sense 
that it satisfies the Cauchy-Riemann equation $$\CR_{\Ju}(u)=du+\Ju(u)\cdot du \cdot j =0$$ and is asymptotically 
cylindrical over the Reeb orbits in $\Gamma^{\pm}$ near the positive/negative punctures, which is equivalent to finiteness 
of the Hofer energy, $E(u)<\infty$ ([BEHWZ],[EGH]). \\ 

It will become crucial in our discussion that, before we can look for an appropriate compactification, we first need to 
divide out the all the natural symmetries of the domain and the target. While the natural $\IR$-action on the target 
$(\IR\times V,\Ju)$ leads to a natural $\IR$-action on every moduli space $\IM^0$, in the case when the sphere carries less 
three punctures we further need to divide out the automorphism group $\Aut(\Si,j)$ of the domain given by 
$$\varphi.(u,j)=(u\circ\varphi,\varphi^*j),\;\varphi\in\Aut(\Si,j),$$ which is trivial when there 
are at least three punctures. Note that in the latter case the holomorphic curve is called domain-stable. Dividing out these 
obvious symmetries we obtain the desired moduli space of \emph{unparametrized} curves, 
$$\IM=\IM(\Gamma^+,\Gamma^-)=\IM^0(\Gamma^+,\Gamma^-)/(\Aut(\Si,j)\times \IR).$$ 

In [BEHWZ] it is shown that the moduli space $\IM(\Gamma^+,\Gamma^-)$ can be compactified by adding nodal curves as in 
Gromov-Witten theory as well as multi-level curves as in Floer theory. \footnote{One additionally needs to fix an absolute homology class $A\in H_2(V)$ 
which we omit here in order to keep notation simple.} In order to define invariants, it remains to prove 
that the compactified moduli space $\CM=\CM(\Gamma^+,\Gamma^-)$ can be equipped with a nice structure, where it will 
turn out that the best we can hope for is a branched-labelled manifold with boundaries and corners.        

\section{Classical Banach space bundle setup}
In this section we shortly recall the main points of the classical approach which is used (in the hope) to prove nice 
results about the compactified moduli space. Instead of looking at the moduli space $\IM$ of unparametrized curves or even the 
compactified moduli space $\CM$ directly, in the classical approach one first tries to prove a nice result about the 
moduli space $\IM^0$ of parametrized curves. \\

Here the well-known idea is to use an infinite-dimensional version of the implicit function theorem to equip $\IM^0$ with 
a nice manifold (or orbifold) structure. In order to be able to apply the implicit function theorem, we need to find an 
infinite-dimensional (Banach) manifold $\BB^0$ of maps as well as a function, more generally, a section $s$ in an appropriate 
infinite-dimensional (Banach space) bundle $\EE^0$ over this manifold, whose zero set is the moduli space of 
parametrized curves, $\IM^0=s^{-1}(0)$. For manifolds of maps and their construction we refer to [Sch] and the foundational 
paper [E] by Eliasson. \\

Following [BM] (see also [F1], [F2] and [W]) an appropriate Banach manifold $\BB^0$ is given 
by the product $$\BB^0 = H^{1,p,\delta}(\Gamma^+,\Gamma^-) \times \IM_{0,s^++s^-},$$ where $\IM_{0,s^++s^-}$ is 
the moduli space of complex structures on the $s^++s^-$ punctured sphere and $H^{1,p,\delta}(\Gamma^+,\Gamma^-) \subset 
H^{1,p}_{\loc}(\Si,\IR\times V)$ is a space of $H^{1,p}$-maps which again satisfy appropriate asymptotic convergence to the 
chosen Reeb orbits in $\Gamma^{\pm}$ near the positive/negative punctures. Note that $\IM_{0,s^++s^-}$ just consists 
of a point when $s^++s^- = 3$, where we define that this continues to hold when $s^++s^-<3$, i.e., when the curve is 
domain-unstable. Furthermore it can be shown that there exists a smooth infinite-dimensional Banach space bundle $\EE^0$ 
over $\BB^0$ where the fibre $$\EE^0_{u,j}=L^{p,\delta}(\Lambda^{0,1}\otimes_{j,\Ju}u^*T(\IR\times V)),\; (u,j)\in\BB^0$$ 
consists of $(0,1)$-forms (with respect to the complex structures $j$ on $T\Si$ and $\Ju$ on $u^*T(\IR\times V)$) with 
values in the pullback bundle $u^*T(\IR\times V)$. For more details we refer to the papers above. \\

It can further be shown that the Cauchy-Riemann operator $\CR_{\Ju}$ defines a smooth section in the Banach space bundle 
$\EE^0\to\BB^0$ by $\CR_{\Ju}(u,j)=du+\Ju(u)\cdot du \cdot j\in\EE^0_{u,j}$ for $(u,j)\in\BB^0$. Furthermore it can be 
shown that, if the asymptotic convergence conditions for $H^{1,p,\delta}(\Gamma^+,\Gamma^-)$ are chosen appropriately, 
the moduli space $\IM^0=\IM^0(\Gamma^+,\Gamma^-)$ agrees with the zero set of the section $\CR_{\Ju}$ in $\EE^0\to\BB^0$. 
In order to show that the moduli space $\IM^0$ can be equipped with a nice manifold structure, in the very same way as 
for the usual implicit function theorem it hence suffices to prove that zero is a regular value, that is, that the Cauchy-Riemann section meets 
the zero section transversally. \\

After equipping the moduli space $\IM^0$ of parametrized curves with a nice manifold structure, we deduce from the 
smoothness of the action of the automorphism group on $\IM^0$ that also the moduli space $\IM$ of unparametrized curves 
carries a nice structure. Finally, it follows from the compactness theorem in the Gromov-Hofer topology proven in [BEHWZ], 
together with a gluing theorem which again relies on the above transversality assumption, that the compactified moduli 
space carries a manifold with boundary and corner structure, which in turn can be used to deduce the desired algebraic 
statements. \\ 
       
\section{Problems with multiply-covered curves}
Up to small modifications it can be shown as for Gromov-Witten theory, see [D], that for a generic choice of cylindrical almost 
complex structure the required transversality assumption is true at every \emph{simple} $\Ju$-holomorphic curve 
$(u,j)\in\IM^0$. In other words, assuming that the moduli space $\IM^0$ consists only of curves $(u,j)$, for which there 
does \emph{not} exist a covering map $\varphi: (\Si,j)\to(\Si',j')$ of punctured Riemann surfaces such that 
$u=v\circ\varphi$, then one can prove that the linearization of the Cauchy-Riemann operator $\CR_{\Ju}$ 
is surjective at every $(u,j)\in\IM^0$ for a generic choice of $\Ju$. In particular, it follows that the moduli space 
$\IM^0$ of parametrized \emph{simple} $\Ju$-holomorphic curves carries a nice manifold structure. \\
   
In order to prove this fundamental transversality result for simple curves one shows that the space of almost complex structures is indeed 
large enough to ensure that the Cauchy-Riemann operator is transversal to the zero section when one additionally allows 
the almost complex structure to vary. Forgetting the (simple) holomorphic curve in the resulting universal moduli space 
but only remembering the almost complex structure, one can prove that for every almost complex structure which is regular 
for this forgetful map the resulting linearized Cauchy-Riemann operator is surjective, so that the result follows from 
Sard's lemma.  \\

While it is not directly clear from this transversality proof for simple $\Ju$-holomorphic curves that 
the transversality result is neccessarily false for multiply-covered curves, it is however not very hard to find examples of multiply-covered 
curves where transversality can definitively not be satisfied. \\

Indeed, let $u=v\circ\varphi:(\Si,j)\to(\IR\times V,\Ju)$ be a multiply-covered curve with underlying simple curve 
$v:(\Si',j')\to(\IR\times V,\Ju)$ and branched covering map $\varphi:(\Si,j)\to(\Si',j')$ and assume that for the 
Fredholm indices for $u$ and $v$ we have $$\ind(u)<\ind(v)+2\cdot\#\Crit(\varphi),$$ where $\Crit(\varphi)$ is the set 
of branch points of the branching map $\varphi$. Since the kernel of the linearized operator $D_u$ contains at least the 
infinitesimal variations of $u$ as a multiply-covered curve, it follows that 
$$\dim\ker D_u \geq \dim\ker D_v +2\cdot\#\Crit(\varphi) \geq \ind(v)+2\cdot\#\Crit(\varphi),$$
and hence $\dim\ker D_u > \ind(u)$, while we need to get equality when transversality would be satisfied. \\

On the other hand, it is not very hard to find examples of multiply-covered curves for which the above inequality for the 
Fredholm indices of the multiple cover and the underlying simple curve is satisfied. Indeed, as shown in [F2], the above 
problem already occurs for the basic examples of $\Ju$-holomorphic curves in $\IR\times V$, namely the branched covers 
of orbit cylinders $\IR\times\gamma$ over closed Reeb orbits. \\

While the underlying simple curve, the orbit cylinder, always 
has Fredholm index zero, it can be shown using standard estimates for the Conley-Zehnder index for multiply-covered orbits 
that for every $\Ju$-holomorphic curve $u$ which is a branched cover of an orbit cylinder, we always have 
$\ind(u)\leq 2\cdot\#\Crit(\varphi)$. While it can be shown that the right hand side agrees with $\dim\ker D_u$ and is 
independent of the underlying closed Reeb orbit, the Fredholm index $\ind(u)$ crucially relies on the underlying orbit, 
that is, the Conley-Zehnder index of $\gamma$ and its iterates. It can be seen that in general we have to expect to 
get a strict inequality, which in turn shows that transversality fails in general even for the basic examples of 
holomorphic curves studied in SFT. \\

As a concrete example, choose a hyperbolic orbit in a three-dimensional contact manifold 
and consider a holomorphic curve with two positive and one negative puncture (pair-of-pants) branching over the 
corresponding orbit cylinder such that it is asymptotically cylindrical over $\gamma$ near the two positive punctures and 
over the double-covered orbit $\gamma^2$ near the negative puncture: while $\dim\ker D_u = 2\cdot\#\Crit(\varphi)=2$, an easy index 
computation yields $\ind(u)=1$. \\ 

Note that these orbit curves actually cause trouble in two ways. First, it can be shown using an area estimate, 
that when the chosen collections of closed Reeb orbits $\Gamma^{\pm}$ only involve iterates 
$\gamma^k$ of the same closed Reeb orbit $\gamma$, the corresponding moduli spaces $\IM(\Gamma^+,\Gamma^-)$ entirely 
consist of orbit curves over the corresponding orbit cylinder. To be more precise, it is shown in [F2] that there 
exists a natural action filtration on SFT such that in the corresponding spectral sequence the first differential only 
counts orbit curves. \\

While it was shown by the author in [F2] (and [F3] after introducing gravitational descendants) how one has to 
deal with these moduli spaces of orbit curves, the orbit curves further cause trouble in another way. Indeed, while 
it follows from the compactness theorem in [BEHWZ] that an appropriate subsequence of a sequence of $\Ju$-holomorphic curves 
in $\IM(\Gamma^+,\Gamma^-)$ converges to a $\Ju$-holomorphic curve with multiple levels and nodes, there exists no 
estimate to exclude that some levels consist of orbit curves. While it is natural to expect this to happen for a 
general moduli space of $\Ju$-holomorphic curves, it is one of the key problems in the embedded contact homology (ECH) of 
Taubes and Hutchings that orbit curves even appear in the compactification of the nice moduli spaces of embedded 
holomorphic curves studied in this subtheory of SFT. \\

Indeed, the compactness statement proven by Hutchings and Taubes (based on [BEHWZ]) shows that a sequence of embedded 
curves with index two converge to a holomorphic curve with many 
levels, where only the first and the last level consist of embedded curves of index one, while all the other levels in 
between consist of orbit curves of index zero. It follows that, in order to prove that $\del\circ\del=0$ in ECH, the 
authors had to prove in [HT1] and [HT2] a generalized gluing formula for holomorphic curves in SFT, where curves are glued after possibly 
inserting additional levels of orbit curves of index zero. \\

While at first glance this seems to contradict the compactness statement used to prove the algebraic statements in 
SFT, the difference just reflects the fact that transversality can not be established for all moduli spaces using 
generic choices of almost complex structures, which is however crucial for positivity of intersections 
needed in the definition of ECH.               

\section{Special approaches to the transversality problem}
After illustrating how severe the transversality problem in SFT actually is, we want to recall in this section three 
ways to achieve transversality in special cases. Note that none of these approaches solves the transversality 
problem in SFT in full generality. On the other hand, we decided not to discuss other approaches like 
virtual fundamental cycle constructions and Kuranishi structures in this note which claim to prove the 
transversality problem for holomorphic curves in full generality.  

\subsection{Automatic transversality in dimension four}
It is well-known that holomorphic curve techniques are particularly powerful for four-dimensional symplectic manifolds 
and, in the same way, for three-dimensional contact manifolds. The reason is that, since in this case the holomorphic 
curves have half the dimension of the surrounding symplectic manifold, one can use intersection theory to prove stronger 
results. Following [MDSa] the central observation is that holomorphic curves always intersect positively with respect to the induced 
orientations, so that the number of intersection points and singularities can be bounded from above by purely 
topological quantities. \\
 
The intersection theory is used in two ways: apart from the fact that it excludes that sequences of 
'sufficiently nice' curves converge to 'bad' multi-level curves, say, with orbit curve components, it implies automatic 
transversality results which even hold for multiply-covered curves and non-generic choices of almost complex structures. \\

To be more precise, based on previous work by Hofer-Lizan-Sikorav, R. Siefrings and V. 
Shevishisin, C. Wendl has proven in [W] an automatic transversality result for SFT in dimension four, which 
he also used to prove that sequences of so-called nicely-embedded curves, which are the building blocks of finite energy 
foliations in SFT in dimension four (see [W] and the reference therein), converge to multi-level curves without orbit 
curve components. In particular, while none of the definitions is stronger than the other, it follows that the embedded 
curves studied by Wendl (which, roughly speaking, have embedded image in the contact manifold $V$) satisfy a stronger compactness result 
than the embedded curves studied in embedded contact homology by Hutchings and Taubes (which have, again roughly speaking, 
embedded image in $\IR\times V$). \\

Apart from the fact that his results only apply when the underlying contact manifold 
is three-dimensional, it only establishes transversality for very special classes of holomorphic curves, since problems 
occur when the curve has genus, singularities or approaches Reeb orbits with even Conley-Zehnder index. On the other hand, 
A. Momin has shown in his recent preprint [M] that, in case there do not exist contractible Reeb orbits, one can define 
cylindrical contact homology for three-dimensional contact manifolds by combining Wendl's automatic transversality result 
with the standard transversality result for simple curves.  

\subsection{Obstruction bundles}
As we have seen in the last section it is easy to find examples of multiply-covered curves where transversality 
can never be satisfied. In other words, it follows from index computations that at such curves the linearized 
Cauchy-Riemann operator must always have a nonzero cokernel. In these cases, it follows that the best one can hope 
for is that the moduli space of multiply-covered curves is a manifold of the wrong dimension and that the cokernels 
fit together to give finite-dimensional bundle over this moduli space, where the rank of the bundle is given by the 
difference between the wrong dimension of the moduli space and the right dimension expected by the Fredholm index. 
Following C. Taubes, see also [MDSa] and [LP], one can then use sections of this \emph{obstruction bundle} to perturb 
the Cauchy-Riemann operator 
in such a way that the perturbed Cauchy-Riemann is transversal to the zero section in the Banach space bundle whenever 
the section of the obstruction bundle is transversal to the zero section and the zero set of the obstruction bundle 
section can be identified with the regular moduli space defined by the perturbed Cauchy-Riemann operator.  \\

Let $\IM'$ denote a moduli space of simple $\Ju$-holomorphic curves $v:(\Si',j')\to(\IR\times V,\Ju)$ and consider 
the corresponding moduli space $\IM$ of multiple covers \\$u=v\circ\varphi:(\Si,j)\to(\Si',j')\to(\IR\times V,\Ju)$, where 
$v:(\Si',j')\to(\IR\times V,\Ju)$ is a simple curve in $\IM'$ and $\varphi:(\Si,j)\to(\Si',j')$ is any branched 
covering map between punctured Riemann surfaces. For generic choice of almost complex structure $\Ju$ it easily 
follows from the transversality result for simple curves that not only $\IM'$ but also the moduli space $\IM$ 
of multiple covers $u=v\circ\varphi$ carries a nice manifold structure of (local) real dimension 
$\ind(v)+2\#\Crit(\varphi)$. To be more precise, for the tangent spaces we have 
$$T_u\IM=T_v\IM'\oplus \IC^{\#\Crit(\varphi)},$$ where the second summand keeps track of the deformations of 
$u=v\circ\varphi$ as a multiple cover. \\

In order to show that the cokernels of the linearized operators $D_u$ fit together to give a nice vector bundle 
of the right rank, it suffices by $\ind(u) = \dim\ker D_u - \dim\coker D_u$ to prove that $\dim T_u\IM = \dim\ker D_u$, 
which by the obvious inclusion $T_u\IM\subset\ker D_u$ is indeed equivalent to $T_u\IM=\ker D_u$. In other words, 
it remains to show that multiple covers are sufficiently isolated in the sense that every infinitesimal deformation 
of a multiple cover as a holomorphic curve is still a multiply-covered curve, which however 
need not be true and is in general very hard to establish. \\

However, in [F2] it is shown that this is true for 
moduli space of orbit curves, the basic examples of holomorphic curves in SFT, using an infinitesimal version of 
the statement in [BEHWZ] that every curve with zero contact area is an orbit curve. Apart from the fact that in [F2] 
(and [F3] with gravitational descendants) this obstruction bundle method is employed to actually compute the contribution 
of orbit curves to the SFT invariants, in an ongoing project with C. Wendl the author is using the above automatic transversality 
result in dimension four to show that the obstruction bundle technique can be used to study multiple covers of 
nicely-embedded curves. Together with the invariants for closed Reeb orbits obtained in [F2], [F3] this is used to 
define a local version of SFT, where the objects are closed Reeb orbits instead of contact manifolds, while the 
nicely-embedded curves replace the symplectic cobordisms as morphisms of the topological quantum field theory.     

\subsection{Domain-dependent almost complex structures}
The last special approach to the transversality problem we want to discuss briefly are domain-dependent almost complex 
structures. While this is the special approach which is most promising to prove transversality in SFT in full generality, 
so far it could only be successfully applied to the transversality problem in Gromov-Witten theory by K. Cieliebak and 
K. Mohnke in [CM] and to SFT in the Floer case in [F1], while the general case currently seems out 
of reach. \\ 

Roughly speaking, the main idea of this approach is to 'correct' the transversality proof for simple curves in such 
a way that it also proves transversality for multiply-covered curves. Recall that for the transversality result for 
simple curves one shows that the space of almost complex structures is indeed large enough to ensure that the 
Cauchy-Riemann operator is transversal to the zero section when one additionally allows the almost complex structure to 
vary. While the latter is in general not true for multiply-covered curves, it can be shown, see [Sch], that 
transversality indeed 
holds also for multiply-covered curves when one enlarges the dimension of the base by allowing the almost complex 
structure to additionally depend on points on the underlying punctured Riemann surface. \\

While this approach can be used to show that all moduli spaces $\IM^0$ of parametrized curves can be transversally 
cut out of the Banach manifold $\BB^0$ for generic choices of domain-dependent almost complex structures, we would be 
able to solve the transversality problem in SFT when we could ensure that the automorphism groups $\Aut(\Si,j)$ still 
act on $\IM^0$ so that we can still define the desired moduli space $\IM$ of unparametrized curves as quotient, and that 
all the choices of domain-dependent almost complex structures can be made coherent so that everything extends nicely 
to the compactified moduli space $\CM$. \\

It follows from the first requirement that the chosen domain-dependent almost 
complex structures must be covariant with respect to the action of the automorphism group. While for domain-stable curves, 
that is, curves with three or more punctures, this automorphism group is trivial, for domain-independent curves like holomorphic spheres, planes or cylinders 
(without additional marked points) one can show that the allowed domain-dependent almost complex structure are no longer 
allowed to depend on points on the domain, i.e., must be almost complex structures in the usual sense. It follows that 
the approach still cannot establish transversality for curves which are not only multiply-covered but also domain-unstable. \\ 

In [CM] the authors solve this problem for Gromov-Witten theory by employing Donaldson's construction of 
symplectic hypersurfaces in closed symplectic manifolds. Indeed, while holomorphic spheres without marked points 
are domain-unstable, they use these hypersurfaces to introduce enough marked points on the spheres so that they become 
domain-stable. In order to guarantee that in the end one counts the same objects as in usual Gromov-Witten theory, it 
is important that the number of marked points is determined a priori by the topological quantity of energy of the holomorphic sphere, since the hypersurface represents a (large multiple of) the homology class Poincare-dual to the cohomology class given by 
the symplectic form. \\

At first sight the situation in SFT seems to be significantly easier, since almost all holomorphic curves 
studied in SFT are already domain-stable since they carry enough positive and negative punctures. However, if not 
explicitly excluded, domain-unstable curves like spheres, planes or cylinders still exist and need to be stabilized 
in order to prove transversality using domain-dependent almost complex structures. \\

In [F1] it is shown that in the Floer case with symplectically aspherical symplectic manifold $M$ one can establish 
transversality using domain-dependent almost complex structures which are symmetric with respect to the natural 
$S^1$-symmetry on the underlying trivial mapping torus $V=S^1\times M$, where the symmetry is used to compute the SFT-invariants for this 
important special case. While holomorphic spheres and planes do not exist, it is crucial for the proof that for 
holomorphic cylinders, which agree with Floer cylinders, transversality can be established with $S^1$-symmetric 
almost complex structures as long as the Hamiltonian is sufficiently small in the $C^2$-norm such that all Floer 
cylinders are indeed Morse trajectories. \\

Apart from the fact that for the general case it is very hard to exclude the 
existence of holomorphic planes, there already exists a problem with holomorphic cylinders in the general case, 
since as in usual Floer homology transversality for holomorphic cylinders can only be established for cylindrical almost complex
structures which additionally depend on the natural $S^1$-coordinate on the cylinder. It follows that one needs to find a 
way to 'coherently' fix special points on all closed Reeb orbits (not only the simple ones) which in turn can be used to fix unique 
$S^1$-coordinates on all holomorphic cylinders, which does not exist to date.
 
\section{A new Fredholm theory}
While these special approaches are very useful for computations, seeing the limitations in all existing approaches the reader 
is hopefully motivated enough to approach the transversality problem for holomorphic curves using a completely new 
Fredholm theory. Apart from the fact that in this new Fredholm theory transversality for the Cauchy-Riemann operator 
can be achieved by simply perturbing it into general position as well-known for finite-dimensional vector bundles, it 
further also allows to prove that all computations using different special approaches lead to the same result as it 
contains all the perturbations for the Cauchy-Riemann operator obtained from the special approaches as a special case.
  
\subsection{Transversality using abstract perturbations}
It is well-known from classical differential topology that in finite-dimensional vector bundles over finite-dimensional 
manifolds every section can be made transversal to the zero section by slightly perturbing it into general position. 
Instead of using special approaches to achieve transversality like domain-dependent almost complex structures, it is 
clear that the most natural solution to the transversality problem would consist of such an abstract transversality 
result using generic perturbations. \\

Recall that in the classical approach one has to show that the moduli space $\IM^0$ of parametrized curves is 
transversally cut out by the Cauchy-Riemann operator as a section in the infinite-dimensional Banach space bundle 
$\EE^0\to\BB^0$. On the other hand, this moduli space can only be compactified in geometric way after dividing out 
the the action of the automorphism group $\Aut(\Si,j)$ of the domain (and the natural $\IR$-action on the target 
$\IR\times V$) to give the moduli space $\IM$ of unparametrized curves, while the compactified moduli space 
$\CM$ is equipped with a nice manifold structure with boundary and corners using a gluing result, which however relies 
on transversality. \\

Using an abstract transversality result for sections in the Banach space bundle $\EE^0\to\BB^0$ as in the 
finite-dimensional case to obtain a perturbed Cauchy-Riemann operator $\CR_{\Ju}^{\nu}=\CR_{\Ju}+\nu$, it would follow that the 
transversality problem is solved as long as it is guaranteed that 
\begin{enumerate}
\item the automorphism group $\Aut(\Si,j)$ of the domain (and the $\IR$-action on the target) still acts on the resulting 
perturbed moduli space $$(\IM^0)^{\nu}=(\CR_{\Ju}^{\nu})^{-1}(0)\subset\BB^0$$ of parametrized curves to define the perturbed 
moduli space $\IM^{\nu}$ of unparametrized curves and 
\item the different abstract perturbations $\nu$ for the different moduli spaces appearing the compactification of the 
moduli space can be chosen coherently so that the perturbed moduli spaces can be glued together to give a new perturbed 
compactified moduli space $\CM^{\bar{\nu}}$. \\
\end{enumerate}            

First, while it is easy to see that (1) is not satisfied for a general abstract perturbation of the Cauchy-Riemann operator in 
$\EE^0\to\BB^0$, recall from the last sectio that it was precisely this symmetry property which caused the problems 
with domain-unstable curves in the transversality approach using domain-dependent almost complex structures. On the other 
hand, assuming for the moment that this first problem could be solved, it is not clear on the level of abstract 
perturbations, which are still sections in infinite-dimensional Banach space bundles, in which sense they "smoothly 
fit together" such that (2) holds, that is, the resulting perturbed moduli spaces can be glued together to give a 
perturbed compactified moduli space with a nice manifold structure with boundaries and corners. \\

While the problems (1) and (2) are much easier to handle for the special transversality approaches like obstruction 
bundles and domain-dependent almost complex structures, it is not clear how to solve them within the classical approach 
using Banach space bundles over Banach manifolds. In the end this is due to the fact that in both special approaches 
the abstract perturbations actually depend on finite-dimensional spaces. Indeed, while in the obstruction bundle and the 
domain-dependent almost complex structures approaches the resulting perturbations of the Cauchy-Riemann operator depend 
on points on the underlying nonregular moduli space of multiple covers or moduli space of punctured Riemann surfaces, 
respectively, the general abstract perturbations using the abstract perturbation approach depend on points on the 
infinite-dimensional Banach manifold. \\

Before discussing the problems appearing with infinite-dimensional manifolds, let us first sketch our dreamland. \\

First, in order to find perturbations for the Cauchy-Riemann operator using an abstract perturbation scheme solving 
(1), i.e., which are in addition symmetric with respect to the action of the automorphism group of the domain (and 
the $\IR$-symmetry of the target), we would hope to find a new "infinite-dimensional bundle" $\EE\to\BB$, which 
contains the Cauchy-Riemann operator as a "smooth" section, such that the moduli space $\IM$ of \emph{unparametrized} curves 
(not $\IM^0$) is the zero set of this Cauchy-Riemann section. Indeed, perturbing the Cauchy-Riemann operator in this 
new bundle setup into general position, the resulting abstract perturbation automatically has the desired property. \\

On the other hand, proceeding in the same way, we can further solve the problem (2) to find abstract perturbations 
which coherently fit together "in a smooth way", if we can construct an even more general "infinite-dimensional bundle" 
$\CE\to\CB$ which again contains the Cauchy-Riemann operator as "smooth" section, but whose zero 
set is now already the compactified moduli space $\CM$. Again, perturbing the Cauchy-Riemann operator in this 
new bundle setup into general position, we find a "smooth" abstract perturbation which restricts to smooth abstract 
perturbations for the noncompact moduli space and the moduli spaces appearing in the boundary such that the resulting 
perturbed moduli spaces can be glued to a perturbed compactified moduli space with the desired manifold structure.        
   
\subsection{Motivation for sc-structures and sc-manifolds}
In the same way as the moduli space $\IM$ was obtained from the moduli space $\IM^0$ by dividing out the symmetries of 
the domain (and the target), $\IM=\IM^0/\Aut(\Si,j)(\times \IR)$, the canonical way to obtain the infinite-dimensional bundle 
$\EE\to\BB$ solving (1) is to make use of the fact that the action extends to the Banach manifold of maps $\BB^0$ and 
also lifts to an action on the Banach space bundle $\EE^0$ and define $$\EE=\EE^0/\Aut \to \BB=\BB^0/\Aut.$$ 
Provided that the symmetry group $\Aut(\Si,j)$ still acts smoothly with respect to the Banach manifold topology on 
$\EE^0\to\BB^0$, it would follow that the new bundle $\EE\to\BB$ is still a Banach space bundle over a Banach manifold 
and (1) could be proved using the same abstract transversality result as used for the classical Banach space bundle 
$\EE^0\to\BB^0$. \\

In order to solve (1) it hence only remains to verify that the action of $\Aut(\Si,j)$ is indeed smooth with respect to 
the classical Banach manifold topology. As toy model for symmetries of the domain, we restrict for the moment 
our attention to classical Morse theory. Here the Banach manifold of maps is the space of paths from one fixed (critical) 
point on the underlying manifold to another fixed (critical) point, where for simplicity we assume that the manifold 
is just the linear space and the paths start and end at $0\in\IR^n$. \\

In this case the Banach manifold of paths is just 
the linear Banach space $H^{k,p}(\IR,\IR^n)$ with $k\geq 1,p\geq 2$ and we consider the action by translations, 
$$ \tau: \IR \times H^{k,p}(\IR,\IR^n) \to H^{k,p}(\IR,\IR^n),\; (s,u)\mapsto u(\cdot + s).$$

In order to check whether this map is smooth, we first need to check that it is even differentiable. Assuming for the 
moment that $u$ is actually an element in the dense subspace $C^{\infty}_0(\IR,\IR^n)\subset H^{k,p}(\IR,\IR^n)$ of 
smooth functions (with appropriate decay to zero near $\pm\infty$), basic calculus gives us the result 
$$D\tau: \IR \times H^{k,p}(\IR,\IR^n) \to H^{k,p}(\IR,\IR^n),\; D\tau(0,u)\cdot (1,\xi) =u'+\xi,\; 
u'=\frac{\del u}{\del s}. $$ 

Since the formula for the differential of the shift map involves 
the first derivative of the map $u\in H^{k,p}(\IR,\IR^n)$ it follows that, in order to define a map to 
$H^{k,p}(\IR,\IR^n)$, \emph{the differential is only defined at maps $u$ in the subspace 
$H^{k+1,p}(\IR,\IR^n)$ !} \\   

In the same way as in the Morse theory example, it follows that the differential of the action of the symmetries of the 
domain is only defined at elements in the dense subset $\BB^0_1\subset\BB^0$ with 
$$\BB^0_1=H^{k+1,p,\delta}(\Gamma^+,\Gamma^-)\times\IM_{0,s}.$$ Proceeding further, it can be shown that the higher 
derivatives can only be computed at points in the dense subsets $\BB^0_{\ell}= H^{k+\ell,p,\delta}(\Gamma^+,\Gamma^-)\times
\IM_{0,s}$, while all derivatives only exist at maps $(u,j)$ in the dense subset 
$\BB^0_{\infty} = C^{\infty}_0(\Gamma^+,\Gamma^-)\times \IM_{0,s}$ of smooth maps, where 
$$C^{\infty}_0(\Gamma^+,\Gamma^-) := \bigcap_{\ell=0}^{\infty} H^{k+\ell,p,\delta}(\Gamma^+,\Gamma^-).$$ 

It follows that in the classical Banach manifold topology on $\BB^0$ the action of the symmetries is \emph{not} 
smooth, \emph{so that the naive approach for equipping the bundle $\EE=\EE^0/\Aut\to\BB=\BB^0/\Aut$ with the 
structure of a smooth infinite-dimensional bundle to solve (1) does not work !} \\

On the other hand, in order to define a new (infinite-dimensional) differential topology on $\BB^0$ (and $\EE^0$), 
it seems natural to equip the underlying Banach spaces $E$ with a scale structure 
$$...\subset E_2\subset E_1\subset E_0 = E$$ as it is well-known in wavelet theory (also called multi-resolution analysis) in 
numerical mathematics. In order to ensure that the
derivative always exists, the tangent bundle $T\BB^0$ to the manifold of maps $\BB^0$ should only sit over the dense subset 
$\BB^0_1$,$$T\BB^0\to\BB^0_1,$$ so that differentiability only needs to be checked at curves with higher regularity. \\

While all this sounds rather distressing, in the very same way it should come as a great surprise that Hofer, Wysocki 
and Zehnder were actually able to build a new (infinite-dimensional) differential topology on 
the afore-mentioned \emph{scale-calculus} (sc-calculus) (see [HWZ] and [H] and the references therein), which is on the one hand 
\begin{itemize}
\item \emph{weak enough} such that the action of the symmetries on the infinite-dimensional manifolds (in the new sense) is 
      smooth, but on the other hand also 
\item \emph{strong enough} such that the zero set of the Cauchy-Riemann operator is, after perturbing the latter into 
      general position, a finite-dimensional manifold in the usual sense. 
\end{itemize}

In order to see that this new scale differential topology really reduces to the usual differential topology in finite 
dimensions, it should be mentioned that the inclusions in the sequences of subspaces $...\subset E_2\subset E_1\subset 
E_0 = E$ are actually required to be compact, so that neccessarily $E=E_0=E_1=E_2=...$ if $E$ is finite-dimensional. 
In particular, recall that for the finite-dimensional moduli spaces there is no problem with smoothness since all the holomorphic 
curves are automatically smooth by elliptic regularity.    

\subsection{Motivation for retracts and polyfolds}
After this short digression on the problem of equipping the bundle $\EE\to\BB$ with a nice smooth structure in order to 
solve (1) and on the main ideas how Hofer and his collaborateurs were able to solve it using their sc-calculus, it remains 
the problem to equip the bundle $\CE\to\CB$ with an even more general smooth structure in order to 
also solve (2) to obtain a general transversality result for the Cauchy-Riemann operator, which would in turn solve 
the transversality problem in SFT in full generality. \\

Assume that a sequence of holomorphic curves $(u_n,j_n)$ in $\IM$ converges to a broken or nodal curve 
$(u_{\infty},j_{\infty})$, where $u_{\infty}: (\Si_{\infty},j_{\infty})\to(\IR\times V,\Ju)$ starts from the 
nodal Riemann surface with punctures $(\Si_{\infty},j_{\infty})$. We want to describe a neighborhood of 
$(u_{\infty},j_{\infty})$ in the space of maps $\CB$ containing the compactified moduli space $\CM$. \\

For this we first need to glue the nodal Riemann surface $(\Si_{\infty},j_{\infty})$ to a smooth Riemann surface $(\Si_r,j_r)$ 
depending on the gluing parameter $r\in [0,\infty)$ (for $r$ sufficiently large). Since the gluing 
should be defined as usual by choosing cylindrical 
coordinates $[0,\infty)\times S^1$ and $(-\infty,0]\times S^1$ near the double-points, removing the half-infinite 
cylinders $[r,\infty)\times S^1$ and $(-\infty,-r)\times S^1$ and gluing the resulting Riemann surface with boundary 
in the obvious way, it follows that we can think of the glued smooth Riemann surface $(\Si_r,j_r)$ as a (non-connected !) 
subset of $(\Si_{\infty},j_{\infty})$. Neglecting the use of cut-off function for smoothness for simplicity, it follows 
that there exists a natural (pre-)gluing map to obtain from $(u_{\infty},j_{\infty})$ a family of maps $(u_r,j_r)$ 
starting from the smooth Riemann surface $\Si_r$ by defining $u_r = u_{\infty}|_{\Si_r}$ as $\Si_r\subset\Si_{\infty}$. \\

Forgetting about the variation of the complex structure on the domain, recall that the tangent space 
at $(u_{\infty},j_{\infty})$ to the manifold of maps starting from the nodal Riemann surface $\Si_{\infty}$ should be given by
$$ H^{k,p,\delta}(u_{\infty}^*T(\IR\times V)) = H^{k,p,\delta}(\Si_{\infty},\IC^n), $$
where the identity follows by choosing a unitary trivialization of the pull-back bundle $u_{\infty}^*T(\IR\times V)$. 
On the other hand, for the (pre-)glued curve $(u_r,j_r)$ the space of infinitesimal deformations as a map starting from the 
glued Riemann surface $\Si_r$ is given by 
$$ H^{k,p,\delta}(u_r^*T(\IR\times V)) = H^{k,p,\delta}(\Si_r,\IC^n), $$
where we explicitly assume that the unitary trivialization of $u_r^*T(\IR\times V)$ is given by the unitary trivialization 
of $u_{\infty}^*T(\IR\times V)$ using that $u_r=u_{\infty}|_{\Si_r}$. \footnote{When one uses cut-off functions 
for the definition of the pre-glued map we need to employ parallel transport in order to identify the unitary trivializations.}\\

Viewing $\Si_r$ as a subset of $\Si_{\infty}$ and using zero extension (we again forget about cut-off functions for sufficient 
regularity for simplicity) it follows that we can naturally view the space of infinitesimal variations 
(the tangent space to the manifold of maps) of each (pre-)glued map $(u_r,j_r)$ as a subset of the space of 
infinitesimal variations of the underlying nodal map $(u_{\infty},j_{\infty})$. Moreover, there exists a natural projection 
$$\pi_r: H^{k,p,\delta}(\Si_{\infty},\IC^n) \to H^{k,p,\delta}(\Si_r,\IC^n), $$ which is just given by restricting (we still 
forget about cut-off functions for sufficient regularity for simplicity).  
In particular, observe that, while the moduli spaces which one has to add in order 
to compactify a moduli space always have a smaller dimension than the original moduli space, at least as expected by the 
Fredholm index, the space of variations as a (not neccessarily holomorphic) map is always larger for the curves in 
the boundary. \\

It follows that the desired manifold of maps $\CB$, containing the compactified moduli space $\CM$ 
as a subset, should locally near the nodal curve $(u_{\infty},j_{\infty})$ be given (forgetting the variation of 
the complex structure for simplicity) by
\begin{eqnarray*}
 O &=& \bigcup_{r\in R} \{r\}\times H^{k,p,\delta}(\Si_r,\IC^n) \\
   &=& \{(r,\pi_r(u)): (r,u)\in R\times H^{k,p,\delta}(\Si_{\infty},\IC^n)\} \\
   &=& \rho(U), 
\end{eqnarray*}
when we choose $$U= R\times H^{k,p,\delta}(\Si_{\infty},\IC^n),\; \rho:U\to U,\; \rho(r,u)=(r,\pi_r(u)), $$
so that, in particular, $\rho\circ\rho=\rho$. \\

But with our reasoning we directly arrived at the definition of a \emph{retract} by Hofer, Wysocki and Zehnder, which 
is the local model of a M-polyfold. Using their scale-calculus which we motivated before, they show in their papers 
that these new geometric objects can be equipped with an infinite-dimensional differentiable topology. In particular, 
it turns out that their scale differential topology is again indeed 
\begin{itemize}
\item \emph{weak enough} to allow dimension jumps, in particular, deal with the fact that the dimension goes up at 
the strata which contain the boundary of the moduli space, but still
\item \emph{strong enough} to detect the different strata, in particular, distinguish between interior and boundary strata. 
\end{itemize}

\subsection{The Fredholm package}
Besides defining a new infinite-dimensional differential topology for sc-manifolds and, more generally, (M-)polyfolds, 
Hofer, Wysocki and Zehnder show in their papers that one can generalize the usual Fredholm package used in the 
classical approach to these new categories of geometric objects. In particular, apart from defining (non-linear) 
Fredholm sections in bundles, they prove an implicit function theorem and an abstract perturbation 
scheme leading to the desired abstract transversality result for the Cauchy-Riemann operator. \\

While in the sc-category the Fredholm property is still very natural and leads to an elegant proof of elliptic regularity 
by making use of the sc-structure, for M-polyfolds one is faced with the problem that their local model is no longer a 
linear space but a retract which can be viewed as a parametrized family of subspaces. Since these subspaces have locally 
varying dimensions, more precisely, are obtained from a fixed infinite-dimensional linear space by removing varying 
infinite-dimensional subspaces, it does no longer suffice to call the linearization Fredholm when the linear operator is 
an isomorphism up to a finite-dimensional kernel and cokernel. Due to this problem with the jumping of dimensions, Hofer and 
his collaborateurs introduced the notion of a "filler" and a "filled section", which solve the problem and lead to a 
generalized definition of Fredholm. \\

Finally, while we have already mentioned and used that the new differential topology is 
weaker than the classical differential topology, this has the draw-back that for the new implicit function theorem (IFT) for 
sc-manifolds and polyfolds it does not only suffice to study the linearization of the section, but one more generally 
needs to study the germ of the section (naturally defined using the sc-structure) near the zero. The resulting 
"germ IFT" however suffices to equip the zero set of the Cauchy-Riemann operator with a nice 
manifold structure with boundary and corners, of course, only after possibly perturbing it into general position using 
the abstract perturbation scheme.


\begin{thebibliography}{100000}

\bibitem[BEHWZ]{BEHWZ} Bourgeois, F., Eliashberg, Y., Hofer, H., Wysocki, K. and Zehnder, E.: {\it Compactness results in 
      symplectic field theory.} Geom. and Top. {\bf 7}, 2003. 

\bibitem[BM]{BM} Bourgeois, F. and  Mohnke, K.: {\it Coherent orientations in symplectic field theory.} Math. Z. {\bf 248}, 2003.

\bibitem[CM]{CM} Cieliebak, K. and K. Mohnke: {\it Symplectic hypersurfaces and transversality for Gromov-Witten theory.} 
      J. Symp. Geom. 5, pp. 281-356, 2007.

\bibitem[D]{D} Dragnev, D.: {\it Fredholm theory and transversality for noncompact pseudoholomorphic maps in symplectizations.}
        Comm. Pure Appl. Math. {\bf 57}(6), 2004.

\bibitem[E]{E} Eliasson, H.: {\it The geometry of manifolds of maps.} J. Diff. Geom. {\bf 1}, 1967.

\bibitem[EGH]{EGH} Eliashberg, Y., A. Givental, H. Hofer: {\it Introduction to symplectic field theory.} GAFA 2000 Visions in 
Mathematics special volume, part II, pp. 560-673, 2000.

\bibitem[F1]{F1} Fabert, O.: {\it Contact homology of Hamiltonian mapping tori.} Comm. Math. Helv. {\bf 85}, pp. 203-241, 2010. 

\bibitem[F2]{F2} Fabert, O.: {\it Obstruction bundles over moduli spaces with boundary and the action filtration in 
symplectic field theory.} ArXiv preprint (0709.3312), 2007. 

\bibitem[F3]{F3} Fabert, O.: {\it Gravitational descendants in symplectic field theory.} ArXiv preprint (0907.0789), 2009.

\bibitem[H]{H} Hofer, H.: {\it Polyfolds and a general Fredholm theory.} ArXiv preprint (0809.3753), 2008.

\bibitem[HT1]{HT1} Hutchings, M. and Taubes, C.: {\it Gluing pseudoholomorphic curves along branched covered cylinders I.} 
      J. Symp. Geom. {\bf 5}, 2007.
 
\bibitem[HT2]{HT2} Hutchings, M. and Taubes, C.: {\it Gluing pseudoholomorphic curves along branched covered cylinders II.} 
      J. Symp. Geom. {\bf 7}, 2009.

\bibitem[HWZ]{HWZ} Hofer, H., Wysocki, K. and Zehnder, E.: {\it A general Fredholm theory I: A splicing-based differential geometry.} 
      J. Eur. Math. Soc. {\bf 9}(4), 2007. 

\bibitem[LP]{LP} Lee, J. and T. Parker: {\it An obstruction bundle relating Gromov-Witten invariants of curves and Kahler 
surfaces.} ArXiv preprint (0909.3610), 2009.

\bibitem[M]{M} Momin, A.: {\it Cylindrical contact homology on complements of Reeb orbits.} ArXiv preprint (1003.3268), 2010.

\bibitem[MDSa]{MDSa} McDuff, D. and D. Salamon: {\it $J$-holomorphic curves and symplectic topology.} AMS Colloquium 
Publications 52, 2004.

\bibitem[Sch]{Sch} Schwarz, M.: {\it Cohomology operations from $S^1$-cobordisms in Floer homology.} Ph.D. thesis, Swiss Federal Inst. of 
        Techn. Zurich, Diss. ETH No. 11182, 1995.

\bibitem[W]{W} Wendl, C.: {\it Automatic transversality and orbifolds of punctured holomorphic curves in dimension 
four.} ArXiv preprint (0802.3842), 2008.

\end{thebibliography}
\end{document}